\documentclass[11pt]{amsart}
\usepackage{amssymb,amsmath,amsthm,amsfonts,
}

\newcommand{\comment}[1]{}

\newtheorem{lem}{Lemma}
\newtheorem{propn}[lem]{Proposition}

\newtheorem{thm}[lem]{Theorem}
\newtheorem*{thmN}{Theorem \textnormal{(Lyall \cite{N1})}}

\theoremstyle{remark}
\newtheorem{rem}[lem]{Remark}

\theoremstyle{definition}

\newcommand{\R}{\mathbf R}
\newcommand{\Rn}{{\mathbf R}^d}
\newcommand{\h}{{\mathbf H}_a}

\newcommand{\p}{\partial}

\newcommand{\vp}{\varphi}
\newcommand{\D}{\delta}

\newcommand{\A}{\alpha}

\newcommand{\B}{\beta}
\newcommand{\lm}{\lambda}

\newcommand{\bfa}{\mathbf a}
\newcommand{\bfb}{\mathbf b}
\newcommand{\bfc}{\mathbf c}
\newcommand{\bfd}{\mathbf d}
\newcommand{\bfe}{\mathbf e}
\newcommand{\bff}{\mathbf f}
\newcommand{\bfg}{\mathbf g}

\newcommand{\Tr}{\text{t}}
\DeclareMathOperator{\rank}{rank}

\begin{document}
\title[Strongly Singular Integral Operators on the Heisenberg Group]{Strongly Singular Integral Operators associated to different quasi-norms on the Heisenberg Group}
\author{Norberto Laghi\hspace{3cm}Neil Lyall}
\thanks{Both authors were partially supported by HARP grants from the European Commission; the first author is currently supported by an EPSRC grant, the second author is partially supported by a NSF FRG grant.}
\address{School of Mathematics, The University of Edinburgh, JCM
 Building, The King's Buildings, Edinburgh EH9 3JZ, United Kingdom}
\email{N.Laghi@ed.ac.uk} 
\address{Department of Mathematics, The University of Georgia, Boyd
  Graduate Studies Research Center, Athens, GA 30602, USA}
\email{lyall@math.uga.edu}

\begin{abstract}
In this article we study the behavior of strongly singular integrals associated to three different, albeit equivalent, quasi-norms on Heisenberg groups; these quasi-norms give rise to phase functions whose mixed Hessians may or may not drop rank along suitable varieties. In the particular case of the Koranyi norm we improve on the arguments in \cite{N1} and obtain sharp $L^2$ estimates for the associated operators.
\end{abstract}
\maketitle

\setlength{\parskip}{5pt}

\section{Introduction
}

The Heisenberg group $\h^n$ is a non-commutative nilpotent Lie group, with underlying manifold $\R^{2n+1}$ equipped the group law
\begin{equation}
(x,t)\cdot (y,s)=(x+y,s+t-2a \ x^\Tr J y)
\end{equation}
where $a$ is a 
nonzero real number and $J$ denotes the standard symplectic matrix on $\R^{2n}$, namely
\[J=\left(\begin{matrix} 0 &I_n \\ -I_n &0 \end{matrix}\right)\]
with inverses given by $(x,t)^{-1}=-(x,t)$.
The \emph{nonisotropic} dilations
\begin{equation}
(x,t)\mapsto
(\D x,\D^2 t).
\end{equation}
are automorphisms of $\h^n$ and as such the homogeneous dimension of this group is $2n+2$.

We will consider here, for different quasi-norms $\rho(x,t)$ on $\h^n$, the class of model (group) convolution operators 
formally given by
\[Tf(x,t)=f*K_{\A,\B}(x,t)\]
where $K_{\A,\B}$ is a strongly singular distributional kernel on $\h^n$ 
that agrees, for $(x,t) \ne (0,0)$, with the function
\[K_{\A,\B}(x,t) = {\rho(x,t)}^{-2n-2-\alpha}e^{i{\rho(x,t)}^{-\beta}}\chi(\rho(x,t)),\]
where 
$\B>0$ and $\chi$ is smooth and compactly supported in a small neighborhood of the origin
.

Operators of this type were first studied in the Euclidean setting of $\R^d$ with Fourier transform techniques (and $\rho(x)=|x|$) by Hirschman \cite{Hi} in the case $d=1$ and then in higher dimensions by Wainger \cite{W}, Fefferman \cite{F1}, and Fefferman and Stein \cite{FS}. For some generalizations and an oscillatory integral approach to these classical results, see Lyall \cite{N2}. 

The analogous questions on the Heisenberg group were first investigated by the second author in \cite{N1} using the group Fourier transform and some `partial' results were obtained. In this article we optimally sharpen the previously obtained results and in addition address the question of the behavior of different quasi-norms for the first time. We do not employ the Fourier transform in our arguments and as such our methods are not restricted to the class of translation-invariant operators.
 
There are various natural choices of quasi-norm on $\h^n$,
for example one can take
\[\rho_0(x,t)=\max\{|x_1|,\dots,|x_{2n}|,|t|^{1/2}\},\]
which is however not smooth away from the origin. We shall instead consider the following three equivalent quasi-norms which clearly are smooth away from the origin;
\begin{itemize}
\item[(i)] $\rho_1(x,t)=(|x|^4+t^2)^{1/4}$ 
\medskip
\item[(ii)] $\rho_2$ defined by $\rho_2(x,t)=1\Longleftrightarrow |x|^2+t^2=1$ and extended by homogeneity
\medskip
\item[(iii)] $\rho_3(x,t)=(x_1^4+\cdots+x_{2n}^4+t^2)^{1/4}$ 
.
\end{itemize}


The case when $\rho(x,t)=\rho_1(x,t)$, the so called Koranyi norm on the Heisenberg group, was initially studied by the second author in \cite{N1} using the group Fourier transform. 
\begin{thmN} Let $\rho(x,t)=\rho_1(x,t)$ and $a\ne0$.
\begin{itemize}
\item[(i)] If $\A\leq n\B$, then $T$ extends to a bounded operator from $L^2(\h^n)$ to itself.
\item[(ii)] If $T$ extends to a bounded operator from $L^2(\h^n)$ to itself, then $\A\leq (n+\frac{1}{2})\B$.
\end{itemize}
\end{thmN}
We note that this result is \emph{uniform} in $a\ne0$. Its proof relied upon the radial nature of the Koranyi norm and uniform asymptotic expansions for Laguerre functions due to Erd\'elyi \cite{E} along with some careful analysis.

The first main result of this article is the following.

\begin{thm}\label{A} There exists a constant $C_\B$ such that if $\rho(x,t)=\rho_1(bx,bt)$ with $0<a^2/b^2<C_\B$, then $T$ extends to a bounded operator from $L^2(\h^n)$ to itself if and only if $\A\leq (n+\frac{1}{2})\B$.
\end{thm}

In actual fact we show that
\[C_\B=\frac{\B+2}{2}\left(2\B+5+\sqrt{(2\B+5)^2-9}\right)\] 
is admissible and note that we consequently have $C_\B\geq9$ for all $\B>0$.

The analogous result for  $\rho(x,t)=\rho_2(x,t)$, the nonisotropic Minkowski functional associated to the (Euclidean) unit ball, is the following.

\begin{thm}\label{B} If $\rho(x,t)=\rho_2(bx,bt)$ with $a^2/b^2\leq1$, then $T$ extends to a bounded operator from $L^2(\h^n)$ to itself whenever $\A\leq (n+\frac{1}{2})\B$. 
\end{thm}

We note that the Koranyi norm $\rho_1(x,t)$ is indeed an actual norm on $\h^n$ for all $0<a^2\leq1$, as is the nonisotropic Minkowski functional $\rho_2(x,t)$ for sufficiently small $|a|$.

A negative result related to Theorems \ref{A} and \ref{B} for 
$\rho_3(x,t)$ is discussed in Section \ref{detcal}.

\section{Reduction to dyadic estimates}

The necessary condition in Theorem \ref{A} follows from the arguments in \cite{N1}.
To establish sufficiency in both Theorems \ref{A} and \ref{B} matters reduce to considering the dyadic operators \[T_j(x,t)=f*K_j(x,t),\] where
\[K_j(x,t)=\vartheta(2^j\rho(x,t))K_{\A,\B}(x,t),\]
with $\vartheta\in C_{0}^{\infty}(\R)$ supported in $[\frac{1}{2},2]$ is chosen such that $\sum_{j=0}^{\infty}\vartheta(2^jr)=1$ for all $0\leq r\leq1$.

As in \cite{N1} everything reduces to establishing the following key dyadic estimates.

\begin{thm}\label{Tj} If $\A\leq (n+\frac{1}{2})\B$ and either 
\[\text{(i) } \rho(x,t)=\rho_1(bx,bt) \text{ with } 0<a^2/b^2<C_\B \quad
\text{or}\quad \text{(ii) } \rho(x,t)=\rho_2(bx,bt) \text{ with } a^2/b^2\leq1\] 
then the dyadic operators $T_j$ are bounded uniformly on $L^2(\h^n)$, more precisely
\begin{equation}\label{L2}
\int_{\h^n}|T_{j}f(x,t)|^2\, dx\, dt \leq C2^{j(2\A-(2n+1)\B)}\int_{\h^n}|f(x,t)|^2\, dx\, dt.
\end{equation}
\end{thm}

Theorems \ref{A} and \ref{B} then follow from an application of Cotlar's lemma (and a standard limiting argument) since the operators $T_j$ are, in the following sense, almost orthogonal.
\begin{propn}\label{T*T}
If $\A\leq(n+\frac{1}{2})\B$, then 
\[\|T^*_jT_{j'}\|_{L^2(\h^n)\rightarrow L^2(\h^n)}+\|T_jT^*_{j'}\|_{L^2(\h^n)\rightarrow L^2(\h^n)}\leq C2^{-(n+\frac{1}{2})\B|j-j'|}.\]
\end{propn}
This follows exactly as in \cite{N1} once we have made the observation that if $\rho(x,t)$ is any quasi-norm on $\h^n$ satisfying the estimate $c^{-1}\leq\rho(\D x,\D^2t)\leq c$ for some $c\geq1$, then there exists a constant $c_0>1$ so that either
\[c_0^{-1}\leq\frac{\p}{\p x_j}\rho(x,t)\leq c_0\]
for some $j=1,\dots,2n$, or
\[c_0^{-1}\D\leq\frac{\p}{\p t}\rho(x,t)\leq c_0\D.\]
Since from this it follows that
\[\left|\nabla_{(y,s)}[\rho(y,s)^{-\B}-\rho((x,t)\cdot (y,s))^{-\B}]\right|\geq C>0,\]
whenever $\rho((x,t)\cdot (y,s))\gg \rho(y,s)$. For more details see \cite{N1}.

\section{Homogeneous groups and a proposition of H\"ormander}

The Heisenberg group is of course one of the simplest examples of a (non-commutative)  homogeneous group.
Recall that a homogeneous group consists of $\Rn$ equipped with a Lie group structure, together with a family of dilations
\[x=(x_1,\dots, x_d)\mapsto\D\circ x=(\D^{a_1}x_1,\dots,\D^{a_d}x_d),\]
with $a_1,\dots,a_d$ strictly positive, that are group automorphisms, for all $\D>0$. 

To each homogeneous group on $\Rn$, we can associate its Lie algebra, consisting of left-invariant vector fields on $\Rn$, with basis $\{X_j\}_{1\leq j\leq d}$ where each $X_j$ is the left-invariant vector field that agrees with $\p/\p x_j$ at the origin.

Key to establishing Theorem \ref{Tj} is the following, presumably well known, generalization of a proposition of H\"ormander \cite{FLp}, see also 
\cite{BigS}, Chapter IX.

\begin{propn}\label{vcp}
Let $\Psi$ be a smooth function of compact support in $x$ and $y$, and $\Phi$ be real-valued and smooth on the support of $\Psi$. If we assume that 
\begin{equation}\label{nd}\det\Bigl(X_j Y_k\Phi(x,y)\Bigr)\neq0,\end{equation}
on the support of $\Psi$, then for $\lm>0$ we have
\begin{equation}\label{H} \Bigl\|\int_{\Rn}\Psi(x,y)e^{i\lm\Phi(x,y)}f(y)dy\Bigr\|_{L^2(\Rn)}\leq C\lm^{-\frac{d}{2}}\|f\|_{L^2(\Rn)}.\end{equation}
\end{propn}

Proposition \ref{vcp} can in fact be extended to families of smooth vector fields $X_1,\dots,X_d$ forming a basis at every point of $\R^d$; however, we choose to state it in this restricted generality (which is already more than we need) as this admits a proof which is simply the natural modification of H\"ormander's original argument.

\begin{proof}
By using a partition of unity we may assume that the amplitude $\Psi$ has suitably small compact support in both $x$ and $y$. Denoting the operator on the left hand side of inequality (\ref{H}) by $T_\lm$
it is then easy to see that
\[T_\lm^*T_\lm f(y)=\int_{\Rn}K_\lm(x,z)f(z)\,dz\]
where
\[K_\lm(x,z)=\int_{\Rn} e^{i\lm[\Phi(x,y)-\Phi(z,y)]}\Psi(x,y)\overline{\Psi(z,y)}\,dy.\]
It therefore suffices to establish the kernel estimate
\begin{equation}\label{Kest}
|K_\lm(x,z)|\leq C(1+\lm|z^{-1}\cdot x|)^{-N},
\end{equation}
since from this it would follow that
\[\int |K_\lm(x,z)|\,dz\approx|\{z:|z^{-1}\cdot x|\leq\lm^{-1}\}|=C\lm^{-d}\]
and similarly for $\int |K_\lm(x,z)|\,dx$, and therefore by Schur's test that 
\[\|T_\lm^*T_\lm f\|_{L^2(\Rn)}\leq C \lm^{-d} \|f\|_{L^2(\Rn)}.\]

The kernel $K_\lm(x,z)$ is of course always bounded, hence in order to establish (\ref{Kest}) we need only consider the case when $|z^{-1}\cdot x|\geq\lm^{-1}$. 
Now
\begin{align*}
Y_k\Phi(x,y)-Y_k\Phi(z,y)&=\int^1_0\frac{d}{dt}Y_k\Phi(z\cdot t(z^{-1}\cdot x),y)\,dt\\
&=\sum_{j=1}^d (z^{-1}\cdot x)_j \int^1_0 X_jY_k\Phi(z\cdot t(z^{-1}\cdot x),y)\,dt\\
&=\sum_{j=1}^d (z^{-1}\cdot x)_j X_jY_k\Phi(x,y) +\mathit{O}(|z^{-1}\cdot x|^2).
\end{align*}
So if we let
\[A=A(x,y)=X_j Y_k\Phi(x,y) \ \ \ \ \ \text{and} \ \ \ \ \ \ u=u(x,y,z)=A^{-1}\frac{z^{-1}\cdot x}{|z^{-1}\cdot x|}\]
and define
\[\Delta(x,y,z)=(u_1Y_1+\cdots+u_dY_d)[\Phi(x,y)-\Phi(z,y)]\]
it follows that
\[\Delta(x,y,z)=|z^{-1}\cdot x|+\mathit{O}(|z^{-1}\cdot x|^2).\]
Therefore for $|z^{-1}\cdot x|$ small enough, it is here that we use our initial suitably small support assumption, we have
\[|\Delta(x,y,z)|\geq \tfrac{1}{2}|z^{-1}\cdot x|,\]
and if we now set \[D=\frac{1}{i\lm\Delta(x,y,z)}(u_1Y_1+\cdots+u_dY_d),\]
it follows that
\begin{align*}
\left|\int_{\Rn} e^{i\lm[\Phi(x,y)-\Phi(z,y)]}\Psi(x,y)\overline{\Psi(z,y)}\,dy\right|&=\left|\int_{\Rn} D^N\left(e^{i\lm[\Phi(x,y)-\Phi(z,y)]}\right)\Psi(x,y)\overline{\Psi(z,y)}\,dy\right|\\
&=\left|\int_{\Rn} e^{i\lm[\Phi(x,y)-\Phi(z,y)]}(D^\Tr)^N\left(\Psi(x,y)\overline{\Psi(z,y)}\right)\,dy\right|\\
&\leq C_N(1+\lm|z^{-1}\cdot x|)^{-N},
\end{align*}
for all $N\geq0$.
\end{proof}

\section{Proof of Theorem \ref{Tj}}\label{dyadic}

Since the operator norms of $T_j$ are equal to that of the rescaled operator $\widetilde{T}_j$, given by
\begin{align*}
\widetilde{T}_jf(x,t)= \int_{\h^n} \widetilde{K}_j\bigl((y,s)^{-1}\cdot(x,t)\bigr)f(y,s)\,dy\,ds
\end{align*}
where
\begin{align*}
\widetilde{K}_j(x,t)&=2^{-j(2n+2)}K_j(2^{-j}x,2^{-2j}t)\\
&=2^{j\A}\vartheta(\rho(x,t))\rho(x,t)^{-2n-2-\alpha}e^{i2^{j\B}\rho(x,t)^{-\beta}}
\end{align*}
it suffices to establish estimate (\ref{L2}) for the rescaled operators $\widetilde{T}_j$. 

\comment{
We have already reduced matters to establishing the estimate
\begin{equation}\label{Sj}
\int_{\h^n}|S_{j}f(x,t)|^2\, dx\, dt \leq C2^{-j(2n+1)\B)}\int_{\h^n}|f(x,t)|^2\, dx\, dt
\end{equation}
for the rescaled operators $S_j$.
}
\comment{
Since the operator norms of $T_j$ are equal to that of the rescaled operator $\widetilde{T}_j$, given by
\begin{align*}
\widetilde{T}_jf(x,t)= \int \widetilde{K}_j\bigl((y,s)^{-1}\cdot(x,t)\bigr)f(y,s)\,dy\,ds
\end{align*}
where
\begin{align*}
\widetilde{K}_j(x,t)&=2^{-jd_h}K_j(2^{-j}x,2^{-2j}t)\\
&=2^{j\A}\vartheta(\rho(x,t))\rho(x,t)^{-d_h-\alpha}e^{i2^{j\B}\rho(x,t)^{-\beta}}
\end{align*}
it suffices to establish estimate (\ref{L2}) for the rescaled operators $\widetilde{T}_j$. 
}

Since the $\widetilde{T}_j$ are local operators, in the sense that the support of $\widetilde{T}_jf$ is always contained in a fixed dilate of some nonisotropic ball containing the support of $f$,
we may make the additional assumption that the integral kernels above have compact support in both $(x,t)$ and $(y,s)$. 
Estimate (\ref{L2}) for $\widetilde{T}_j$ then follows from Proposition \ref{vcp} once we have verified the \emph{non-degeneracy} condition (\ref{nd}) in this setting.

It is well known that
\[X^\ell_j=\frac{\p}{\p x_j}+2ax_{j+n}\frac{\p}{\p t}, \ \ X_{j+n}^\ell=\frac{\p}{\p x_{j+n}}-2ax_{j}\frac{\p}{\p t} \ \ \ \ \ j=1,\dots,n,\]
and $T=\frac{\p}{\p t}$ form a real basis for the Lie algebra of left-invariant vector fields on $\h^n$, while 
\[X^r_j=\frac{\p}{\p x_j}-2ax_{j+n}\frac{\p}{\p t}, \ \ X_{j+n}^r=\frac{\p}{\p x_{j+n}}+2ax_{j}\frac{\p}{\p t} \ \ \ \ \ j=1,\dots,n,\]
and $T=\frac{\p}{\p t}$ form a real basis for the Lie algebra of right-invariant vector fields.

For convenience we shall use synonymously
$X^\ell_{2n+1}=X^r_{2n+1}=T,$
and furthermore denote
\[X^\ell=(X^\ell_1,\dots,X^\ell_{2n+1}) \ \ \ \text{and} \ \ \ \ X^r=(X^r_1,\dots,X^r_{2n+1}).\]
We note that
\[-[X_j^r\widetilde{\vp}](x,t)=[X_j^\ell\vp]\bigl((x,t)^{-1}\bigr),\]
where $\widetilde{\vp}(x)=\vp\bigl((x,t)^{-1}\bigr)$, and hence
\[X^\ell_jY^\ell_k\left[\vp\bigl((y,s)^{-1}\cdot (x,t)\bigr)\right]=-[X^\ell_jX^r_k\vp]\bigl((y,s)^{-1}\cdot (x,t)\bigr).\]
The \emph{non-degeneracy} condition (\ref{nd}) in this setting is therefore equivalent to the following.

\begin{propn}\label{KonHn} Let $\Phi(x,t)=\rho(x,t)^{-\beta}$ with $\B>0$. If $(x,t)\ne(0,0)$ and either \[\text{(i) } \rho(x,t)=\rho_1(x,t) \text{ with } 0<a^2<C_\B \quad
\text{or}\quad \text{(ii) } \rho(x,t)=\rho_2(x,t) \text{ with } a^2\leq1\] 
then
\[\det\Bigl(X^\ell_j X^r_k\Phi(x,t)\Bigr)\neq0.\]
\end{propn}
\comment{
\begin{propn}\label{SWonHn} Let $\Phi(x,t)=\rho_2(x,t)^{-\beta}$, then for all $\B>0$ and $(x,t)\ne(0,0)$, 
\[\det\Bigl(X^\ell_j X^r_k\Phi(x,t)\Bigr)\neq0\]
provided $a^2\leq1$.
\end{propn}
}
Theorem \ref{Tj}  now follows immediately for $b=1$, the proof in general follows from the observation that $\h^n$ is isomorphic to $\mathbf{H}^n_{a/b}$ with the explicit isomorphism being given by \[\phi(x,t)=(bx,bt).\]

\section{The determinant calculations}\label{detcal}

The purpose of this section is to prove Proposition \ref{KonHn}, however we shall start by stating and sketching the proof of a related negative result for the quasi-norm $\rho_3(x,t)$ on $\mathbf{H}^1_a$. Outlining this argument first will be instructive as it is simpler than, while still similar to, those for Proposition \ref{KonHn}.

\begin{propn}\label{Bad} Let $n=1$ and $\Phi(x,t)=\rho_3(x,t)^{-\beta}$, then 
\[\det\Bigl(X^\ell_j X^r_k\Phi(x,t)\Bigr)=0\]
along the lines $(0,x_2,0)$ and $(x_1,0,0)$.
\end{propn}
\begin{proof}
Let $\vp_3(x,t)=\rho_3(x,t)^4=x_1^4+x_2^4+t^2$.
It is straightforward to see that the `mixed' Hessian of $\Phi$ is given by
\[X^\ell_j X^r_k\Phi(x,t)=-\tfrac{\B}{4}\vp_3^{-(\B+8)/4}\{\vp_3 X^\ell_j X^r_k\vp_3-\tfrac{\B+4}{4} X^\ell_j\vp_3 X^r_k\vp_3\}.\]

For convenience both here and in the proofs of both parts of Proposition \ref{KonHn} we define
\[A:=X^\ell_jX^r_k\vp_3 \ \ \ \ \ \text{and} \ \ \ \ \ \ B:=X^\ell_j\vp_3 X^r_k\vp_3.\]
Now since $\rank(B)=1$ it follows that
\[\det(\vp_3 A-\tfrac{\B+4}{4}B)=\vp_3^{2}\left\{\vp_3\det(A)-\frac{\B+4}{4}\left\{\det\left(\begin{matrix} \bfb_{1}\\ \bfa_{2}\\ \bfa_{3}\end{matrix}\right)+\det\left(\begin{matrix} \bfa_{1}\\ \bfb_{2}\\ \bfa_{3}\end{matrix}\right)+\det\left(\begin{matrix} \bfa_{1}\\ \bfa_{2}\\ \bfb_{3}\end{matrix}\right)\right\}\right\},\]
where $\bfa_j=(a_{j1},a_{j2},a_{j3})$ and  $\bfb_j=(b_{j1},b_{j2},b_{j3})$.

It is easy to verify that
\[X_1^\ell\vp_3(x,t)=4(x_1^3+ax_2x_3)\quad X_2^\ell\vp_3(x,t)=4(x_2^3-ax_1x_3)\quad\text{and}\quad X_3^\ell\vp_3(x,t)=2t,\]
while
\[X_1^r\vp_3(x,t)=4(x_1^3-ax_2x_3)\quad X_2^r\vp_3(x,t)=4(x_2^3+ax_1x_3)\quad\text{and}\quad X_3^r\vp_3(x,t)=2t,\]
and hence that
\[A=2(C+aD)\]
where
\[C=\left(\begin{matrix} 6x_1^2 & 2t & 0 \\ -2t & 6x_2^2 & 0 \\ 0 & 0& 0\end{matrix}\right)\quad\text{and}\quad D=\left(\begin{matrix} -4ax_2^2 & 4ax_1x_2 & 2x_2 \\ 4ax_1x_2 & -4ax_1^2 & -2x_1 \\ -2x_2 & 2x_1& 1/a\end{matrix}\right)\]

Since $\rank(D)=1$ it follows that
\[\det(A)=8\det\left(\begin{matrix} 6x_1^2 & 2t \\ -2t & 6x_2^2\end{matrix}\right)=32(9x_1^2x_2^2+t^2).\]
For the first of the remaining three determinants we note that
\[\det\left(\begin{matrix} \bfb_{1}\\ \bfa_{2}\\ \bfa_{3}\end{matrix}\right)=8X_1^\ell\vp_3\det(E+aF),\]
where
\[E=\left(\begin{matrix} 2x_1^3 & 2x_2^3 &0\\ -2t & 6x_2^2 & 0 \\ 0&0&0\end{matrix}\right)\quad\text{and}\quad F=\left(\begin{matrix} -2x_2t & 2x_1t &t/a\\ 4ax_1x_2 & -4ax_1^2 & -2x_1 \\ -2x_2&2x_1&1/a\end{matrix}\right).\]
Using the fact that $\rank(F)=1$ we then see that
\[\det\left(\begin{matrix} \bfb_{1}\\ \bfa_{2}\\ \bfa_{3}\end{matrix}\right)=8X_1^\ell\vp_3\det\left(\begin{matrix} 2x_1^3 & 2x_2^3 \\ -2t & 6x_2^2\end{matrix}\right).\]
In an almost identical manner we can also obtain that 
\[\det\left(\begin{matrix} \bfa_{1}\\ \bfb_{2}\\ \bfa_{3}\end{matrix}\right)=8X_2^\ell\vp_3\det\left(\begin{matrix} 6x_1^2 &2t\\2x_1^3 & 2x_2^3 \end{matrix}\right),\]
and
\[\det\left(\begin{matrix} \bfa_{1}\\ \bfa_{2}\\ \bfb_{3}\end{matrix}\right)=8X_3^\ell\vp_3\left\{2ax_2\det\left(\begin{matrix} -2t & 6x_2^2\\2x_1^3 & 2x_2^3 \end{matrix}\right)+2ax_1\det\left(\begin{matrix} 6x_1^2 &2t\\2x_1^3 & 2x_2^3 \end{matrix}\right)+t\det\left(\begin{matrix} 6x_1^2 & 2t \\ -2t & 6x_2^2\end{matrix}\right)\right\},\]
we leave the details to the reader. Bringing all of this together we get that
\[\det(\vp_3 A-\tfrac{\B+4}{4}B)=-16\vp_3^{2}\left\{6(\B+1)\vp_3x_1^2x_2^2+(\B+2)t^4+3(\B+4)x_1^2x_2^2t^2-2(x_1^4+x_2^4)t^2\right\}.\qedhere\]
\end{proof}

\subsection{Proof of Proposition \ref{KonHn}, part (i)} 
Let $\vp_1(x,t)=\rho_1(x,t)^4=|x|^4+t^2$.
It is straightforward to see that the `mixed' Hessian of $\Phi$ is given by
\[X^\ell_j X^r_k\Phi(x,t)=-\tfrac{\B}{4}\vp_1^{-(\B+8)/4}\{\vp_1 X^\ell_j X^r_k\vp_1-\tfrac{\B+4}{4} X^\ell_j\vp_1 X^r_k\vp_1\}.\]

We again define
$A:=X^\ell_jX^r_k\vp_1$ and $B:=X^\ell_j\vp_1 X^r_k\vp_1$.
Since $\rank(B)=1$ it follows that
\[\det(\vp_1 A-\tfrac{\B+4}{4}B)=\vp_1^{2n}\left\{\vp_1\det(A)-\frac{\B+4}{4}\sum_{j=1}^{2n+1}\det\left(\begin{matrix} \bfa_{1}\\ \vdots \\ \bfb_{j}\\
\vdots\\ \bfa_{2n+1}\end{matrix}\right)\right\},\]
where $\bfa_j=(a_{j1},\dots,a_{j\,2n+1})$ and  $\bfb_j=(b_{j1},\dots,b_{j\,2n+1})$.

It is an easy calculation to see that
\[X^\ell\vp_1(x,t)=\bigl(4|x|^2x+4at(Jx),2t\bigr),\]
\[X^r\vp_1(x,t)=\bigl(4|x|^2x-4at(Jx),2t\bigr),\]
where $J$ is the standard symplectic matrix on $\R^{2n}$ coming from the group structure. Hence we have
\[A=4\left(\begin{matrix} C &0 \\ 0 &0 \end{matrix}\right)+8\left(\begin{matrix} D &0 \\ 0 &0 \end{matrix}\right)+4aE,\ \ \ \text{and} \ \ \ \ B=4|x|^2\left(\begin{matrix} F &0 \end{matrix}\right)+4atG,\]
where
\[C=|x|^2I+atJ \ \ \ \ \ D=xx^\Tr \ \ \ \ \ \ E=\left(\begin{matrix} 2a(Jx)(x^\Tr J) &Jx \\ x^\Tr J &1/2a\end{matrix}\right),\]
\[F=(X^\ell\vp_1)x^\Tr \ \ \ \ \ \text{and} \ \ \ \ \ G=\left(\begin{matrix} (X^\ell\vp_1)(x^\Tr J) &X^\ell\vp_1/2a \end{matrix}\right).\]

Now since both $\rank(D)=1$ and $\rank(E)=1$ it follows that
\begin{align*}
\det(A)&=2^{4n+1}\det(C+2D)\\&=2^{4n+1}\left\{\bigl(|x|^4+a^2t^2\bigr)^n+\tfrac{1}{2}\bigl(|x|^4+a^2t^2\bigr)^{n-1}\sum_{j=1}^{2n}x_jX^\ell_j\vp_1\right\}\\
&=2^{4n+1}\bigl(|x|^4+a^2t^2\bigr)^{n-1}\bigl(3|x|^4+a^2t^2\bigr).
\end{align*}
To obtain the final identity above we used the fact that 
\[\sum_{j=1}^{2n}x_jX^\ell_j\vp_1=4|x|^4.\]

We note that for all $j$ for which it makes sense, both
\begin{equation}\label{fandg}
\rank\left(\begin{matrix} \bfe_{1}\\ \vdots \\ \bfg_{j}\\
\vdots\\ \bfe_{2n+1}\end{matrix}\right)=1 \quad \text{and}\quad \rank\left(\begin{matrix} \bfd_{1}\\ \vdots \\ \bff_{j}\\
\vdots\\ \bfd_{2n}\end{matrix}\right)=1.\end{equation} 
From this and the observation that
\[\sum_{j=1}^{2n}(X^\ell_j\vp_1)^2=16|x|^2\bigl(|x|^4+a^2t^2\bigr)\]
we may therefore conclude that
\begin{align*}
\sum_{j=1}^{2n}\det\left(\begin{matrix} \bfa_{1}\\ \vdots \\ \bfb_{j}\\
\vdots\\ \bfa_{2n+1}\end{matrix}\right)&=2^{4n+1}\sum_{j=1}^{2n}\det\left\{\left(\begin{matrix} \bfc_{1}\\ \vdots \\ \mathbf{0}_{j}\\
\vdots\\ \bfc_{2n}\end{matrix}\right)+\left(\begin{matrix} 2\bfd_{1}\\ \vdots \\ |x|^2\bff_{j}\\
\vdots\\ 2\bfd_{2n}\end{matrix}\right)\right\}\\
&=2^{4n+1}|x|^2\sum_{j=1}^{2n}X^\ell_j\vp_1\det\left(\begin{matrix} \bfc_{1}\\ \vdots \\ \mathbf{0}_j+x^\Tr\\
\vdots\\ \bfc_{2n}\end{matrix}\right)\\
&=2^{4n-1}\bigl(|x|^4+a^2t^2\bigr)^{n-1}|x|^2 \sum_{j=1}^{2n}(X^\ell_j\vp_1)^2\\
&=2^{4n+3}|x|^4\bigl(|x|^4+a^2t^2\bigr)^{n}.
\end{align*}

Finally, we can combine (\ref{fandg}) and the fact that
\[\det(C+2D)= \bigl(|x|^4+a^2t^2\bigr)^{n-1}\bigl(3|x|^4+a^2t^2\bigr),\]
together with the identity
\[\sum_{j=1}^{n}(x_{j}X^\ell_{j+n}\vp_1-x_{j+n}X^\ell_j\vp_1)=-4a|x|^2t\]
to obtain
\begin{align*}
\det&\left(\begin{matrix} \bfa_{1}\\ \vdots \\ \bfa_{2n}\\ \bfb_{2n+1}\end{matrix}\right)
=4^{2n+1}a\, X^\ell_{2n+1}\vp_1 \det\left(\begin{matrix} C+2D &Jx\\ |x|^2 x^\Tr & t/2a\end{matrix}\right)\\
&=4^{2n+1}t \left\{2a|x|^2 \sum_{j=1}^{n}\left\{x_{j}\det\left(\!\!\begin{matrix} \bfc_{1}\\ \vdots \\ \mathbf{0}_{j+n}+x^\Tr \\ \vdots \\ \bfc_{2n} \end{matrix}\!\!\right)-x_{j+n}\det\left(\!\!\begin{matrix} \bfc_{1}\\ \vdots \\ \mathbf{0}_{j}+x^\Tr \\ \vdots \\ \bfc_{2n} \end{matrix}\!\!\right)\right\}+t\det(C+2D)\right\}\\
&=4^{2n+1}t\bigl(|x|^4+a^2t^2\bigr)^{n-1} \left\{\tfrac{1}{2}a|x|^2 \sum_{j=1}^{n}(x_{j}X^\ell_{j+n}\vp_1-x_{j+n}X^\ell_j\vp_1)+t\bigl(3|x|^4+a^2t^2\bigr)\right\}\\
&=4^{2n+1}t^2 \bigl(|x|^4+a^2t^2\bigr)^{n-1}\left((3-2a^2)|x|^4+a^2t^2\right).
\end{align*}

Bringing this all together we see that
\[\sum_{j=1}^{2n+1}\det\left(\begin{matrix} \bfa_{1}\\ \vdots \\ \bfb_{j}\\
\vdots\\ \bfa_{2n+1}\end{matrix}\right)=4^{2n+1}\bigl(|x|^4+a^2t^2\bigr)^{n-1}\left\{2|x|^8+t^2\bigl(3|x|^4+a^2t^2\bigr)\right\},\]
and consequently
\[\det(\vp_1 A-\tfrac{\B+4}{4}B)
=-(4\vp_1)^{2n}\bigl(|x|^4+a^2t^2\bigr)^{n-1}f_1(x,t)\]
where
\begin{align*}
f_1(x,t)&=2(\B+1)|x|^8+(\B+2)t^2\bigl(3|x|^4+a^2t^2\bigr)-2|x|^4a^2t^2\\
&=2(\B+1)|x|^8+\bigl(3(\B+2)-2a^2\bigr)|x|^4 t^2+(\B+2)a^2t^4.
\end{align*}
By analyzing the discriminant 
\[\Delta=4a^4-4(\B+2)(2\B+5)a^2+9(\B+2)^2,\]
we see that our Hessian will be \emph{non-degenerate} provided either
\[2a^2\leq3(\B+2) \ \ \ \ \ \text{or} \ \ \ \ \ \ |2a^2-(2\B+5)(\B+2)|<(\B+2)\sqrt{(2\B+5)^2-9},\]
which reduces simply to the condition that
\[a^2<C_\B=\frac{\B+2}{2}\left(2\B+5+\sqrt{(2\B+5)^2-9}\right).\]

\begin{rem}
We conclude by remarking that when $a^2\geq C_\B$ the Hessian degenerates along the paraboloids
\[|x|^4=\frac{2a^2-3(\B+2)\pm\sqrt{\Delta}}{4(\B+1)}\,t^2.\]
In particular when $a^2=C_\B$ we have that $\Delta=0$ and hence the Hessian degenerates along the paraboloid
\[|x|^4=\frac{2C_\B-3(\B+2)}{4(\B+1)}\,t^2=\frac{(\B+1)(\B+2)+\sqrt{(\B+1)(\B+4)}}{2(\B+1)}\,t^2.\]
\end{rem}

\subsection{Proof of Proposition \ref{KonHn}, part (ii)}

We start by letting $\vp_2(x,t)=\rho_2(x,t)^2$ and noting that as a consequence $\vp_2$ must satisy the identity
\begin{equation}\label{def}
\vp_2^{-1}|x|^2+\vp_2^{-2}t^2=1.
\end{equation}

Arguing as in the proof of Proposition \ref{KonHn}, part (i) (and using the same notation) we see that $\vp_2$ satisfies
\comment{
It is straightforward to see that the elements of the `mixed' Hessian of $\Phi$ are given by
\[X^\ell_j X^r_k\Phi(x,t)=-\tfrac{\B}{2}\vp_2^{-(\B+4)/2}\{\vp_2 X^\ell_j X^r_k\vp_2-\tfrac{\B+2}{2} X^\ell_j\vp_2 X^r_k\vp_2\}.\]
For convenience we now define
\[A:=X^\ell_jX^r_k\vp_2 \ \ \ \ \ \text{and} \ \ \ \ \ \ B:=X^\ell_j\vp_2 X^r_k\vp_2.\]
Since $\rank(B)=1$ it follows that
}
\[\det(\vp_2 A-\tfrac{\B+2}{2}B)=\vp_2^{2n}\left\{\vp_2\det(A)-\frac{\B+2}{2}\sum_{j=1}^{2n+1}\det\left(\begin{matrix} \bfa_{1}\\ \vdots \\ \bfb_{j}\\
\vdots\\ \bfa_{2n+1}\end{matrix}\right)\right\}.\]

It is an easy calculation to see that
\[\mathcal{A}X^\ell\vp_2(x,t)=\bigl(2\vp_2^{-1}x+4at(Jx),2\vp_2^{-2}t\bigr),\]
\[\mathcal{A}X^r\vp_2(x,t)=\bigl(2\vp_2^{-1}x-4at(Jx),2\vp_2^{-2}t\bigr),\]
where $J$ is the standard symplectic matrix on $\R^{2n}$ coming from the group structure and
\begin{equation}
\mathcal{A}=\vp_2^{-2}|x|^2+2\vp_2^{-3}t^2=\vp_2^{-1}+\vp_2^{-3}t^2.
\end{equation}
A further (somewhat lengthy) calculation then gives that
\[\mathcal{A}^2A=2\mathcal{A}\vp_2^{-2}\left(\begin{matrix} C &0 \\ 0 &0 \end{matrix}\right)+2\mathcal{A}\vp_2^{-4}D,
\]
where
\[C=\vp_2 I+2atJ\quad\text{and}\quad D=[tX^\ell\vp_2-\vp_2(2aJx,1)][tX^r\vp_2-\vp_2(-2aJx,1)]^\Tr.\]

Now since $\rank(D)=1$ it follows that
\begin{equation}
\det(\mathcal{A}^2A)=2^{2n+1}\mathcal{A}^{2n-1}\vp_2^{-(4n+6)}(\vp_2^2+4a^2t^2)^n|x|^4,
\end{equation}
here we used the fact that 
\[\det(C)=(\vp_2^2+4a^2t^2)^n\quad\text{and}\quad \mathcal{A}(tX^\ell_{2n+1}\vp_2-\vp_2)=\mathcal{A}(tX^r_{2n+1}\vp_2-\vp_2)=-\vp_2^{-1}|x|^2.\]

Using the fact that for all $j=1,\dots,2n+1$
\begin{equation}\label{D}\rank\left(\begin{matrix} \bfd_{1}\\ \vdots \\ \mathbf{0}_{j}\\
\vdots\\ \bfd_{2n+1}\end{matrix}\right)=1\quad\text{}\end{equation}
together with the observation that
\begin{equation}\label{C}\det\left(\begin{matrix} \bfc_{1}\\ \vdots \\ \mathbf{0}_j+2\vp_2^{-1}x^\Tr\\
\vdots\\ \bfc_{2n}\end{matrix}\right)=\vp_2^{-1}(\vp_2^2+4a^2t^2)^{n-1}\mathcal{A}X_j^\ell\vp_2\end{equation}
and some simple reductions we may conclude that
\begin{align*}
\mathcal{A}^{4n+2}\sum_{j=1}^{2n}\det\left(\begin{matrix} \bfa_{1}\\ \vdots \\ \bfb_{j}\\
\vdots\\ \bfa_{2n+1}\end{matrix}\right)
&=4^{n}\mathcal{A}^{2n}\vp_2^{-(4n+2)}|x|^2\sum_{j=1}^{2n}X^\ell_j\vp_2\det\left(\begin{matrix} \bfc_{1}\\ \vdots \\ \mathbf{0}_j+2\vp_2^{-1}x^\Tr\\
\vdots\\ \bfc_{2n}\end{matrix}\right)\\
&=4^{n}\mathcal{A}^{2n-1}\vp_2^{-(4n+1)}(\vp_2^2+4a^2t^2)^{n-1}|x|^2\sum_{j=1}^{2n}(\mathcal{A}X_j^\ell\vp_2)^2\\
&=4^{n+1}\mathcal{A}^{2n-1}\vp_2^{-(4n+5)}(\vp_2^2+4a^2t^2)^{n}|x|^4.
\end{align*}
To obtain the final identity above we used the fact that
\begin{equation}\label{sqr}\mathcal{A}^2\sum_{j=1}^{2n}(X^\ell_{j}\vp_2)^2=4\vp_2^{-4}(\vp_2^2+4a^2t^2)|x|^2.\end{equation}

Using fact (\ref{D}) one more time we see that
\[\mathcal{A}^{4n+2}\det\left(\begin{matrix} \bfa_{1}\\ \vdots \\ \bfa_{2n}\\ \bfb_{2n+1}\end{matrix}\right)
=4^{n}\mathcal{A}^{2n+2}\vp_2^{-4n}X^\ell_{2n+1}\vp_2 X^r_{2n+1}\vp_2 \det(C)
+\sum_{j=1}^{2n} \det \left(\begin{matrix} \bfc_{1}\\ \vdots \\ \bfd_j \\ \vdots \\ \bfc_{2n}\\ \bfb_{2n+1}\end{matrix}\right)\]
It follows from (\ref{C}) that for $1\leq j\leq n$ we have 
\begin{align*}\det \left(\begin{matrix} \bfc_{1}\\ \vdots \\ \bfd_j \\ \vdots \\ \bfc_{2n}\\ \bfb_{2n+1}\end{matrix}\right)&=4^n\mathcal{A}^{2n+1}\vp_2^{-(4n+1)}X^\ell_{2n+1}\vp_2(tX^\ell_j\vp_2-2a\vp_2 x_{j+n})\det\left(\!\!\begin{matrix} \bfc_{1}\\ \vdots \\ \mathbf{0}_{j}+2\vp_2^{-1}x^\Tr \\ \vdots \\ \bfc_{2n} \end{matrix}\!\!\right)\\
&=4^n\mathcal{A}^{2n+2}\vp_2^{-4n}(\vp_2^2+4a^2t^2)^{n-1}X^\ell_{2n+1}\vp_2(tX^\ell_j\vp_2-2a\vp_2 x_{j+n})X^\ell_{j}\vp_2
\end{align*}
and similarly
\[\det \left(\begin{matrix} \bfc_{1}\\ \vdots \\ \bfd_{j+n} \\ \vdots \\ \bfc_{2n}\\ \bfb_{2n+1}\end{matrix}\right)=4^n\mathcal{A}^{2n+2}\vp_2^{-4n}(\vp_2^2+4a^2t^2)^{n-1}X^\ell_{2n+1}\vp_2(tX^\ell_{j+n}\vp_2+2a\vp_2 x_{j})X^\ell_{j+n}\vp_2.\]
It then follows from the identity
\[\mathcal{A}\sum_{j=1}^{n}\left\{x_{j}X^\ell_{j+n}\vp_2-x_{j+n})X^\ell_{j}\vp_2\right\}=-8a^2\vp_2^{-1}t|x|^2\]
together with (\ref{sqr}) that
\begin{align*}
\sum_{j=1}^{2n} \det \left(\begin{matrix} \bfc_{1}\\ \vdots \\ \bfd_j \\ \vdots \\ \bfc_{2n}\\ \bfb_{2n+1}\end{matrix}\right)
&=4^{n+1}\mathcal{A}^{2n}\vp_2^{-(4n+4)}(\vp_2^2+4a^2t^2)^{n-1}X^\ell_{2n+1}\vp_2\left\{(\vp_2^2+4a^2t^2)-2\mathcal{A}a^2\vp_2^{3}\right\}|x|^2t.
\end{align*}

Bringing this all together we see that
\begin{align*}\mathcal{A}^{4n+2}\sum_{j=1}^{2n+1}\det\left(\begin{matrix} \bfa_{1}\\ \vdots \\ \bfb_{j}\\
\vdots\\ \bfa_{2n+1}\end{matrix}\right)&=4^{n+1}\mathcal{A}^{2n}\vp_2^{-(4n+4)}(\vp_2^2+4a^2t^2)^n t^2+\heartsuit
\end{align*}
where
\begin{align*}
\heartsuit&=4^{n+1}\mathcal{A}^{2n-1}\vp_2^{-(4n+5)}(\vp_2^2+4a^2t^2)^{n-1}|x|^2\left\{(\vp_2^2+4a^2t^2)(|x|^2+2\vp_2^{-1}t^2)-4\mathcal{A}a^2\vp_2^2t^2\right\}\\
&=4^{n+1}\mathcal{A}^{2n}\vp_2^{-(4n+1)}(\vp_2^2+4a^2t^2)^{n-1}|x|^2.
\end{align*}
In the last step we have used the identity 
\[\mathcal{A}\vp_2^2=|x|^2+2\vp_2^{-1}t^2.\]

We therefore have that
\[\mathcal{A}^{4n+2}\det(\vp_2 A-\tfrac{\B+2}{2}B)=-2^{2n+1}\mathcal{A}^{2n-1}\vp_2^{-(2n+5)}(\vp_2^2+4a^2t^2)^{n-1}f_2(x,t)\]
where
\begin{align*}
f_2(x,t)&=2\mathcal{A}\vp_2^4|x|^2-(\vp_2^2+4a^2t^2)(|x|^4-2\mathcal{A}\vp_2 t^2)\\
&=\mathcal{A}\vp_2^4|x|^2+4\vp_2^2t^2(1-a^2)+16a^2t^4+4a^2\vp_2^{-2}t^6.
\end{align*}


\section{Further comparisons}

\subsection{Nonisotropic $\R^{2n+1}$}
When $a=0$ (not a Heisenberg type group, but still a homogeneous group) we of course have $X^\ell_j=X^r_j=\p/\p x_j$ and it is then straightforward to verify that in this case we have the following: 

If $\rho(x,t)=\rho_1(x,t)$, then (in the notation of Section \ref{detcal}) we get that
\[\det(\vp_1 A-\tfrac{\B+4}{4}B)=-b\,(4\vp_1)^{2n}|x|^{4n}\left\{2(\B+1)|x|^4+3(\B+2)t^2\right\}\]
and in particular the Hessian degenerates along the line $(0,t)$. 

While if $\rho(x,t)=\rho_2(x,t)$, then
\[\mathcal{A}^{4n+2}\det(\vp_2 A-\tfrac{\B+2}{2}B)=-2^{2n+1}\mathcal{A}^{2n-1}\vp_2^{-5}\left\{\B\mathcal{A}\vp_2^3+\mathcal{A}\vp_2^2|x|^2+4t^2\right\}\]
which is clearly non-degenerate. For a closed related result, see Shayya \cite{BS}.

\subsection{Polarized Heisenberg group}

The polarized Heisenberg group $\mathbf{H}^n_{a,\text{pol}}$ is isomorphic to the full Heisenberg group $\h^n$ and has the multiplication law
\[(x,t)\cdot (y,s)=(x+y,s+t-2a \ x^\Tr J_{\text{pol}} \ y)\]
where again $a$ is a 
nonzero real number, but now $J_{\text{pol}}$ denotes the matrix on $\R^{2n}$, 
\[J_{\text{pol}}=\left(\begin{matrix} 0 &I_n \\ 0 &0 \end{matrix}\right).\]

In particular, if $n=1$ and $a=-1/2$, then this is the $m=3$ case of the groups of $m\times m$ upper-triangular matrices with ones along the diagonal; see Stein \cite{BigS}.

When $n=1$ the corresponding Lie algebra is generated by the left-invariant vector fields 
\[X^\ell_1=\frac{\p}{\p x_1}, \ \ X_{2}^\ell=\frac{\p}{\p x_2}+2ax_{1}\frac{\p}{\p t}, \ \text{ and }X_{3}^\ell=\frac{\p}{\p t}, \]
with the right-invariant vector fields being given by
\[X^r_1=\frac{\p}{\p x_1}+2ax_{1}\frac{\p}{\p t}, \ \ X_{2}^r=\frac{\p}{\p x_2}, \ \text{ and }X_{3}^r=\frac{\p}{\p t}. \]

\begin{propn}\label{BadPol} Let $n=1$ and $\Phi(x,t)=\rho_1(x,t)^{-\beta}$, then for all real $a$
\[\det\Bigl(X^\ell_j X^r_k\Phi(x,t)\Bigr)=0\]
along the line $(0,t)$.
\end{propn}

In actual fact, using again the notation of the previous section, we have that
\[\det(\vp_1 A-\tfrac{\B+4}{4}B)=-16\left\{2(\B+1)|x|^8+3(\B+2)|x|^4t^2-2(\B+2)a\vp_1x_1x_2t \right\}.\]
We leave the details to the reader.

\begin{propn}\label{GoodPol} Let $n=1$ and $\Phi(x,t)=\rho_2(x,t)^{-\beta}$ with $\B>0$, then for $(x,t)\ne(0,0)$
\[\det\Bigl(X^\ell_j X^r_k\Phi(x,t)\Bigr)\ne0\]
provided $a^2\leq1$.
\end{propn}
\begin{proof}
A calculation similar to those above yields
\[\mathcal{A}^{6}\det(\vp_2 A-\tfrac{\B+2}{2}B)=-8\mathcal{A}\vp_2^{-7}\left\{\B\mathcal{A}\vp_2(\vp_2^2-ax_1x_2t)+(2\mathcal{A}\vp_2^3-2\mathcal{A}\vp_2ax_1x_2t-(x_1^2+x_2^2)^2)\right\}.\]
Now it is easy to see that
\begin{align*}
2\mathcal{A}\vp_2^3-2\mathcal{A}\vp_2ax_1x_2t-(x_1^2+x_2^2)^2&=2\vp_2^2+2t^2-(\vp_2-\vp_2^{-1}t^2)^2-2\mathcal{A}\vp_2ax_1x_2t\\
&\geq \vp_2^2-\vp_2^{-2}t^4-2\mathcal{A}\vp_2ax_1x_2t\\
&=\mathcal{A}\vp_2(\vp_2(x_1^2+x_2^2)-2ax_1x_2t)\\
&\geq \mathcal{A}\vp_2(x_1^2+x_2^2)(\vp_2-|at|),
\end{align*}
which is clearly nonnegative if $a^2\leq1$, since it follows from (\ref{def}) that $\vp_2\geq |t|$. In the last line above we used the easily verifiable inequality 
\begin{equation}
2ax_1x_2t\leq(x_1^2+x_2^2)|at|.
\end{equation}
From this inequality it also follows that
\[\vp_2^2-2ax_1x_2t\geq2|x_1x_2|(\vp_2-|at|)\]
and hence, if we again assume that $a^2\leq1$, we see that
\[\vp_2^2-ax_1x_2t\geq\vp_2^2/2.\qedhere\]
\end{proof}

\bibliographystyle{siam}


\end{document}